\documentclass[11pt]{amsart}
\usepackage{epsfig}
\usepackage{graphics}

\newtheorem{theorem}{Theorem}
\newtheorem{lemma}[theorem]{Lemma}
\newtheorem{proposition}[theorem]{Proposition}
\newtheorem{corollary}[theorem]{Corollary}

\theoremstyle{definition}

\theoremstyle{remark}

 \def\Z{{\mathbb{Z}}}

\def\mod{{\rm Mod}}

\def\htc{\mathcal{HT}}

\begin{document}
\newenvironment{prooff}{\medskip \par \noindent {\it Proof}\ }{\hfill
$\square$ \medskip \par}
    \def\sqr#1#2{{\vcenter{\hrule height.#2pt
        \hbox{\vrule width.#2pt height#1pt \kern#1pt
            \vrule width.#2pt}\hrule height.#2pt}}}
    \def\square{\mathchoice\sqr67\sqr67\sqr{2.1}6\sqr{1.5}6}
\def\pf#1{\medskip \par \noindent {\it #1.}\ }
\def\endpf{\hfill $\square$ \medskip \par}
\def\demo#1{\medskip \par \noindent {\it #1.}\ }
\def\enddemo{\medskip \par}

\def\qed{~\hfill$\square$}

\title[Automorphisms of the Hatcher-Thurston complex]
{Automorphisms of the Hatcher-Thurston complex}
\author{Elmas Irmak}
\address{Department of Mathematics, University of Michigan,
Ann Arbor, MI 48109, USA} \email{eirmak@umich.edu}

\author{Mustafa Korkmaz}
\address{Department of Mathematics, Middle East Technical University,
06531 Ankara, Turkey} \email{korkmaz@arf.math.metu.edu.tr}


\date{\today}

\thanks{The first author is supported by a Rackham Faculty Fellowship,
Horace H. Rackham School of
Graduate Studies, University of Michigan. The second author is supported
in part
by the Turkish Academy of Sciences
under the Young Scientists Award Program (MK/T\"UBA-GEB\.IP 2003-10).}

\begin{abstract} Let $S$ be a compact, connected, orientable surface of
positive genus. Let $\htc(S)$ be the Hatcher-Thurston complex of
$S$. We prove that ${\rm Aut}\,\htc(S)$ is isomorphic to the
extended mapping class group of $S$ modulo its center.
\end{abstract}

\maketitle \setcounter{secnumdepth}{1} \setcounter{section}{0}

\section{Introduction}
Let $S$ be a compact, connected, orientable surface of genus $g$
with $r \geq 0$ boundary components. The extended mapping class
group, $\mod_S^*$, of $S$ is the group of isotopy classes of all
homeomorphisms (including orientation reversing) of $S$.
The group $\mod_S^*$ can be viewed as the automorphism group of various
geometric objects. These objects include the complex of curves,
the complex of nonseparating curves, the complex of
separating curves, the complex of pants decompositions and the
complex of Torelli geometry.

The Hatcher-Thurston complex $\htc(S)$, which is defined in Section~\ref{section2} below,
plays a special role in the theory of mapping class groups.
This complex was constructed in~\cite{Ht} by A. Hatcher and W. Thurston
(and used by B. Wajnryb~\cite{w1})
in order to find a presentation for the mapping class group. It was also
used by J. Harer~\cite{Ha} in his computation of the second homology group of
mapping class group. There is a natural action of $\mod_S^*$ on $\htc(S)$
by automorphisms. The purpose of this paper is to show that every automorphism of $\htc(S)$
is induced by some element of $\mod_S^*$. More precisely,
we prove that the automorphism group of the Hatcher-Thurston
complex $\htc(S)$ is isomorphic to the group $\mod_S^*$ modulo its center.
We do this by proving that the automorphism group of our
complex is isomorphic to the automorphism group of the complex $\mathcal{G}(S)$ on
nonseparating simple closed curves, the complex defined by
P. Schmutz Schaller in~\cite{Sc}. (See Section~\ref{section2} for the definition of
$\mathcal{G}(S)$.)

Another complex, the so-called complex of curves $C(S)$, was introduced about
the same time by W. Harvey~\cite{H}. It was also proved to be of
the fundamental importance in the topology of surfaces and in
the theory of Teichm\"uller spaces. Its automorphisms were investigated
in the pioneering paper of N. Ivanov~\cite{Iv}, who proved that the group
of automorphisms of $C(S)$ is equal to the extended mapping class group
of $S$ (for genus~$>1$), and found important applications of this result
to the mapping class groups and to the Teichm\"uller spaces. His result
was used to find automorphisms groups of various other objects related to
surfaces (see \cite{bm}, \cite{cc},\cite{fi}, \cite{Ir3}, \cite{M}, \cite{MW}, \cite{Sc})
and inspired some generalizations
(see \cite{Ir1}, \cite{Ir2}, \cite{K}, \cite{L}).
The paper of D. Margalit~\cite{M} deals with an object closest to the one
considered by us, namely with the so-called pants complex. He proved that
the automorphism group of the pants complex is isomorphic to
the extended mapping class group.
While the results of~\cite{M} and of this paper are similar in the spirit,
neither of them implies the other.

Here is how we prove our main result, Theorem~\ref{thmmain}. The vertices
of the Hatcher-Thurston complex $\htc(S)$ are cut systems.
We encode nonseparating simple closed
curves by vertices and edges of $\htc(S)$. A great deal of work is on this choice.
Using this coding, we show that the automorphism group of $\htc(S)$ has a
well-defined action on the set of (isotopy classes of) nonseparating simple
closed curves. Under this action, we show that dual circles are mapped to
dual circles, giving rise to a homomorphism from the group of automorphisms of $\htc(S)$ to
that of the complex $\mathcal{G}(S)$. We prove that this homomorphism is
in fact an isomorphism. We want to point out that our proof and the proof
of the main result of~\cite{M} are independent of each other, but have some similarities.
(The similarities are pointed out to us by Margalit.)
The similarities are perhaps not surprising, as these methods are indeed very natural
to use in this situation.

\indent The paper is organized as follows. In
Section~\ref{section2}, we give the definition of various
complexes used in the paper and state the relevant properties of
these complexes. In Section~\ref{section3}, we show that every
automorphism $f$ of the Hatcher-Thurston complex $\htc(S)$ induces
an automorphism $\tilde{f}$ of the complex $\mathcal{G}(S)$.
Finally, in Section~\ref{section4}, we state and prove the main
theorem, and discuss some alternative approach.

\section{Various complexes on curves}
\label{section2}

A simple closed curve on $S$ is said to be \textit{nontrivial} (or
\textit{essential}) if it does not bound a disk on $S$ and it is
not homotopic to a boundary component of $S$. We denote simple
closed curves by capital letters and their isotopy classes by the
corresponding lowercase letters. The \textit{geometric intersection number}
$i(a,b)$ of two classes $a$ and $b$ is defined as the minimum number
of intersection points of $A$ and $B$ for $A \in a$ and $B\in b$.

We denote by $\mathcal{A}$ the set of isotopy classes of
nontrivial simple closed curves on $S$. If $C$ is a simple closed
curve on $S$, the surface obtained from $S$ by cutting along $C$
is denoted by $S_C$. Any two simple closed curves $A, B$ are
always assumed to intersect each other minimally. We say that two
simple closed curves $A$ and $B$ on $S$ are \textit{dual} if
they intersect each other transversely at only one point. In this
case we also say that their isotopy classes $a$ and $b$ are dual.
We note that for a simple closed curve $A$ there is a curve dual to $A$
if and only if $A$ is nonseparating.

\subsection{The Hatcher-Thurston complex}

Let $C_1,C_2,\ldots,C_g$ be pairwise disjoint
nonseparating simple closed curves on $S$ such that the surface
obtained from $S$ by cutting along all $C_i$ is connected, so that
it is a sphere with $2g+r$ boundary components.
We call the set $\{ c_1,c_2, \ldots, c_g \}$ a \textit{cut system} and
denote it by $\langle c_1,c_2, \ldots, c_g \rangle$.

Let $v$ and $w$ be two cut systems. Suppose that there are $c\in v$ and
$d\in w$ such that $i(c,d)=1$ and $v-\{ c \}= w-\{ d \}$.
We say that $w$ is obtained from $v$ by an elementary move and we write
$v\leftrightarrow w$.

If $\langle c_1,c_2, \ldots, c_i,\ldots, c_g \rangle
\leftrightarrow \langle c_1,c_2, \ldots, c'_i,\ldots, c_g \rangle$
is an elementary move, then we drop the unchanged curves from the
notation and write $\langle c_i \rangle \leftrightarrow \langle
c'_i \rangle$.

Let $\htc^1(S)$ be the graph obtained by taking cut systems on $S$
as the vertex set and pairs of vertices $\{ v, w\}$ such that
$v\leftrightarrow w$ as the (unordered) edges. This
will be the $1$-skeleton of the Hatcher-Thurston complex.

A sequence of cut systems $(v_1,\ldots, v_n)$ forms a path in
$\htc^1(S)$ if every consecutive pair in the sequence is connected
by an edge in $\htc^1(S)$. A path is closed if $v_1 = v_n$. There
are three types of distinguished closed paths in the graph
$\htc^1(S)$.

\subsubsection{ Triangles} If three vertices have $g-1$ common elements and if the remaining classes
$c,c',c''$ satisfy $i(c,c')=i(c,c'')=i(c',c'')=1$, then
\begin{center}
\unitlength=1mm
\begin{picture}(22,16)
\put(0,12)   {$\langle c \rangle$} \put(10,12) {$\longrightarrow$}
\put(10,12) {$\longleftarrow$} \put(22,12) {$\langle c' \rangle $}
\put(5,5) {$\nwarrow$} \put(5,5) {$\searrow$} \put(18,5)
{$\swarrow$} \put(18,5) {$\nearrow$} \put(10.2,0) {$\langle
c''\rangle$}
\end{picture}
\end{center}
is a triangle (c.f. Figure~\ref{fig1}~(i)). We denote this triangle by
$\langle c \rangle \leftrightarrow \langle c' \rangle \leftrightarrow
\langle c'' \rangle \leftrightarrow \langle c \rangle $.

\subsubsection{ Rectangles}If four vertices have $g-2$ common elements and if the remaining classes
$c_1,c_2,d_1,d_2$ have representatives $C_1,C_2,D_1,D_2$ as in Figure~\ref{fig1}~(ii),
then
\begin{center}
\unitlength=1mm
\begin{picture}(30,16)
\put(0,12)   {$\langle c_1,d_1 \rangle$} \put(25,12)   {$\langle c_1,d_2 \rangle$}
\put(0,0)   {$\langle c_2,d_1 \rangle$} \put(25,0)   {$\langle c_2,d_2 \rangle$}

\put(15,12) {$\longrightarrow$} \put(15,12) {$\longleftarrow$}
\put(15,0) {$\longrightarrow$}  \put(15,0) {$\longleftarrow$}
\put(5,5) {$\updownarrow$}      \put(30,5) {$\updownarrow$}
\end{picture}
\end{center}
is a rectangle. We denote this rectangle by
$\langle c_1,d_1 \rangle \leftrightarrow \langle c_1,d_2 \rangle \leftrightarrow
\langle c_2,d_2 \rangle \leftrightarrow \langle c_2,d_1 \rangle \leftrightarrow
\langle c_1,d_1 \rangle $.

\subsubsection{ Pentagons} 
\begin{figure}[htb]
\begin{center}
\epsfxsize=12cm
\epsfbox{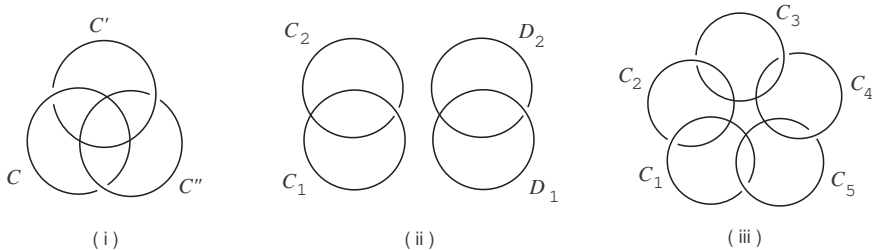}
\caption{A triangle, a rectangle and a pentagon in the Hatcher-Thurston complex}
\label{fig1}
\end{center}
\end{figure}

If five vertices have $g-2$ common elements and if the remaining classes
$c_1,c_2,c_3,c_4,c_5$ have representatives $C_1,C_2,C_3,C_4,C_5$ intersecting each other as in
Figure~\ref{fig1}~(iii), then
\begin{center}
\unitlength=1mm
\begin{picture}(40,34)
\put(7,2) {$\langle c_1,c_4 \rangle$} \put(30,2)   {$\langle
c_1,c_3 \rangle$} \put(0,17)   {$\langle c_2,c_4 \rangle$}
\put(38,17)   {$\langle c_3,c_5 \rangle$} \put(18,30)   {$\langle
c_2,c_5 \rangle$}

\put(21,2) {$\longrightarrow$} \put(21,2) {$\longleftarrow$}
\put(8,9) {$\nwarrow$} \put(8,9) {$\searrow$} \put(38,9)
{$\nearrow$} \put(38,9) {$\swarrow$} \put(10,25) {$\nearrow$}
\put(10,25) {$\swarrow$} \put(35,25) {$\nwarrow$} \put(35,25)
{$\searrow$}

\end{picture}
\end{center}
is a pentagon. Similar to triangles and rectangles we denote this
pentagon by $\langle c_1,c_4 \rangle \leftrightarrow \langle
c_2,c_4 \rangle \leftrightarrow\langle c_2,c_5 \rangle
\leftrightarrow \langle c_3,c_5 \rangle \leftrightarrow\langle
c_1,c_3 \rangle  \leftrightarrow \langle c_1,c_4 \rangle$.

The \textit{Hatcher-Thurston complex} $\htc(S)$ is a two-dimensional
CW-complex obtained from $\htc^1(S)$ by attaching a $2$-cell along each triangle,
rectangle and pentagon.

Hatcher and Thurston used this complex to get a presentation for
the mapping class group for closed orientable surfaces, \cite{Ht}.
They proved that $\htc(S)$ is connected and simply connected.
Wajnryb used it to get a
simple presentation for the mapping class group~\cite{w1} and he
also gave an elementary proof of the connectivity and the
simple connectivity of this complex in~\cite{w2}.\\

\begin{theorem}
\label{thm:htw}
$($\cite{Ht,w2}$)$ Let $S$ be a compact, connected, orientable
surface of genus at least one. Then the complex $\htc(S)$ is
connected.
\end{theorem}

\subsection{The complexes of curves}
The \textit{complex of curves}, $\mathcal{C}(S)$, on $S$ is an
abstract simplicial complex, introduced by Harvey \cite{H}, with
vertex set $\mathcal{A}$, the set of isotopy classes of nontrivial simple closed
curves, such that a set of $n+1$ vertices $\{
a_{0}, a_{1}, a_{2}, \ldots, a_{n} \}$ forms an $n$-simplex
if and only if $a_0, a_1, a_2,\ldots,a_n$ have pairwise disjoint representatives.
The automorphism group of the complex of curves is isomorphic to the extended
mapping class group modulo the center, except for the cases
$(g,r)\in \{ (0,2),(0,3),(0,4),(1,0),(1,1),(1,2)\}$. The reader
is referred to~\cite{Iv},~\cite{K} and~\cite{L} for proof of these results.

Let $\mathcal{B}$ denote the set of isotopy classes of
nonseparating simple closed curves on $S$. The \textit{complex of
nonseparating curves}, $\mathcal{N}(S)$, is the subcomplex of
$\mathcal{C}(S)$ with the vertex set $\mathcal{B}$ such that a set
of $n+1$ vertices $\{ b_0, b_1, b_2,\ldots,b_n \}$ forms an
$n$-simplex if and only if it is an $n$-simplex in
$\mathcal{C}(S)$. If $g \geq 2$, the automorphism group of
$\mathcal{N}(S)$ is isomorphic to the extended mapping class group
of $S$ modulo its center, by the results given in \cite{Ir3}.

In~\cite{Sc}, Schmutz Schaller defined a graph $\mathcal{G}(S)$;
the vertex set is again $\mathcal{B}$, the set of isotopy classes
of nonseparating simple closed curves, and two vertices $a$ and
$b$ are connected by an edge if and only if $i(a,b)=1$. He defines
the graph $\mathcal{G}(S)$ for surfaces of genus zero as well, but
we will not mention that case here. His main result is the following
theorem; we state as much as we need in this paper. Notice that
in the case $g=1$, since the vertices of $\mathcal{G}(S)$ can be viewed
as vertices in $\htc(S)$, the complex $\mathcal{G}(S)$ can be considered
as a subcomplex of $\htc(S)$ in a natural way. In fact, $\mathcal{G}(S)$ is
the $1$-skeleton of $\htc(S)$.

\begin{theorem}
\label{thmsc}
$($\cite{Sc}$)$
Let $S$ be a compact, connected, orientable surface of positive genus.
Then ${\rm Aut}\, \mathcal{G}(S)$ is isomorphic to the extended mapping class group
$\mod^*_S$ modulo the center.
\end{theorem}

In~\cite{Sc}, the case $(g,r)=(1,0)$ is not included, but clearly
it follows from the case $(g,r)=(1,1)$.\\

\subsection{The complex $X_C$}
For a nonseparating simple closed curve
$C$ on $S$, we define a simplicial complex (graph) $X_C$ as
follows: the vertices of $X_C$ are isotopy classes of
nonseparating simple closed curves which are dual to $C$ on $S$. A
set $\{a,b\}$ of vertices forms an edge if and only if $a$ is dual
to $b$.\\

We will need the following definiton in Lemma \ref{lemmaxv}: An
embedded arc $\epsilon$ on a surface $S$ with boundary is called \textit{properly embedded} if
$\partial \epsilon \subseteq \partial S$ and $\epsilon$ is
transversal to $\partial S$. It is called \textit{nontrivial} (or
\textit{essential}) if $\epsilon$ cannot be deformed into
$\partial S$ in such a way that the endpoints of $\epsilon$ stay
in $\partial S$ during the deformation.

\begin{lemma}
\label{lemmaxv} If $S$ is a connected orientable surface of positive genus
and if $C$ is a nonseparating simple closed curve on $S$,
then the complex $X_C$ is connected.
\end{lemma}

\begin{proof} Let $d$ and $d'$ be two distinct vertices in $X_C$.
We will show that there is a path $d = d_0 \rightarrow d_1 \rightarrow
\cdots \rightarrow d_{n+1}=d'$ in $X_C$.

Let $D$ and $D'$ be representatives of $d$ and $d'$ respectively
such that $D$ and $D'$ have minimal intersection and that they are both
dual to $C$. We may assume, moreover, that they intersect $C$ at different points.

If $|D \cap D'|=0$, then $d \rightarrow t_c(d) \rightarrow d'$,
where $t_c$ is the Dehn twist about $c$, is the path in $X_C$ that
we want.

If $|D \cap D'|=1$, then  $d \rightarrow d'$ is the required
sequence.

Assume that $|D \cap D'|= m > 1$. Let $N$ be a regular
neighborhood of $C$ such that the intersection of $D\cup D'$ and
$N$ is a pair of disjoint arcs. Let $R$ be the complement of the
interior of $N$ in $S$ and let $\partial_1$ and $\partial_2$ be
the boundary components of $N$. Let $\epsilon$ and $\tau$ denote
the part of $D$ and $D'$ on $R$ respectively, which are essential
properly embedded arcs. We orient $\epsilon$ and $\tau$ so that
they both start on $\partial_1$ and end on $\partial_2$. We define
an arc in the following way: Start on the boundary component
$\partial_1$ of $R$, on one side of the beginning point of $\tau$
and continue along $\tau$ without intersecting $\tau$, till the
last intersection point of $\epsilon$ and $\tau$ along $\epsilon$.
Then we would like to follow $\epsilon$, without intersecting
$\epsilon \cup \tau$, until we reach $\partial_2$. So, if we are
on the correct side of $\tau$ we do this; if not, we change our
starting side from the beginning and follow the construction. This
gives us an arc, say $\tau_1$. We see that $\tau_1$ is an
essential properly embedded arc since it connects two boundary
components $\partial_1$ and $\partial_2$, and $|\epsilon \cap
\tau_1| < m$ since we eliminated at least one intersection with
$\epsilon$. We also have $|\tau_1 \cap \tau| = 0$ since we never
intersected $\tau$.

Now, using $\epsilon$ and $\tau_1$ in the place of $\epsilon$ and
$\tau$ we define a new properly embedded arc $\tau_2$ connecting
$\partial_1$ to $\partial_2$ such that $|\epsilon \cap \tau_2| <
|\epsilon \cap \tau_1|, |\tau_1 \cap \tau_2|=0$. By an inductive
argument, we get a sequence
$$\epsilon =\tau_{n+1} \rightarrow
\tau_n \rightarrow \tau_{n-1} \rightarrow \cdots \rightarrow \tau_1
\rightarrow \tau_0 = \tau$$
of essential properly embedded arcs on $R$ such that every
consecutive pair is disjoint. So, $\epsilon =\tau_{n+1}$ and
$\tau_n$ are two disjoint arcs on $R$. Note that $\epsilon =
\tau_{n+1}$ is an arc of $D$ on $R$. We now connect the end points
of $\tau_{n}$ with an arc in the interior of $N$ to get a
nonseparating simple closed curve, $Q_{n}$, dual to both $C$ and
$D$ on $S$. Then, we connect the end points of $\tau_{n-1}$ with
an arc in the interior of $N$ so that we get a nonseparating
simple closed curve, $Q_{n-1}$, which is dual to both $C$ and
$\tau_n$ on $S$. Now by an inductive argument, we see that there
is a sequence
$$D = Q_{n+1}\rightarrow Q_n \rightarrow Q_{n-1} \rightarrow \cdots \rightarrow
Q_1 \rightarrow Q_0 = Q$$
consisting of nonseparating simple closed curves dual to $C$ on $S$ such that
every consecutive pair is dual. So, $d$ is connected to $q$ by a
path in  $X_C$.

Since the parts of $D'$ and $Q$ on $R$ are equal, we see that $q =
t_c^m(d')$ for some $m \in \mathbb{Z}$. Then $q$ and $d'$ can be
connected by the path, $q = t_c^m(d') \rightarrow t_c^{m-1}(d')
\rightarrow \cdots \rightarrow t_c(d') \rightarrow d'$ in $X_C$.
Since $d$ and $q$ are connected by a path in $X_C$, we see
that $d$ and $d'$ are connected by a path in $X_C$. Hence, the
complex $X_C$ is connected.\end{proof}

\section{Action of automorphisms of $\htc(S)$ on nonseparating curves}
\label{section3}

We define an action of the automorphism group Aut\,$\htc(S)$
of $\htc(S)$ on the set of nonseparating simple
closed curves as follows. Let $f:\htc (S) \to \htc(S)$ be an
automorphism of the Hatcher-Thurston complex of the surface $S$.
For an isotopy class $c$ of a nonseparating simple closed curve
$C$, choose pairwise disjoint nonseparating simple closed curves
$C_2,C_3,\ldots,C_g$ on $S$ such that $v=\langle c,c_2, \ldots,
c_g \rangle$ is a cut system. Choose another curve $D$ on $S$ such that
$C$ and $D$ are dual and $D$
does not intersect any of $C_i$. Then $w=\langle d, c_2, \ldots,
c_g \rangle$ is also a cut system and the vertices $v$ and $w$ are
connected by an edge in the complex $\htc (S)$. Since $f$ is an
automorphism, the vertices $f(v)$ and $f(w)$ are connected by an
edge as well. Thus the difference $f(v)-f(w)$ of the sets $f(v)$
and $f(w)$ contains only one curve. We define $\tilde{f}(c)$ to be
this unique class.

Notice that if $g=1$ then the cut system $v$ contains only one element; $v=\langle c\rangle $.
Thus $\tilde {f} (c)$ is the unique class in $f\langle c\rangle $, so that
we have $\langle \tilde {f} (c)\rangle = f\langle c\rangle $.

\begin{lemma} \label{lemma1}
For a fixed set of curves $\{ C_2,C_3, \ldots, C_g  \}$, the
definition of $\tilde{f}(c)$ is independent of the choice of the
curve $D$.\end{lemma}

\begin{proof}
Let $v_1\leftrightarrow v_2\leftrightarrow v_3 \leftrightarrow v_1$
be a triangle in the complex $\htc(S)$. Then
we observe that $v_1-v_2=v_1-v_3$.

For a nonseparating simple closed curve $A$ such that
$\langle a,c_2,\ldots,c_g \rangle$ is a cut system,
let $\langle a \rangle$ denote the cut system $\langle a,c_2,\ldots,c_g \rangle$.

Let $D'$ be a simple closed curve on $S$ such that it is dual to both $C$ and $D$, and
disjoint from $C_i$ for $i\geq 2$. Then $\langle d' \rangle =\langle d',c_2,\ldots,c_g \rangle$ is
also a cut system and $ \langle c \rangle \leftrightarrow \langle d
\rangle \leftrightarrow \langle d' \rangle \leftrightarrow \langle
c \rangle$ is a triangle in $\htc(S)$. Since $f$ is an
automorphism, $f\langle c \rangle \leftrightarrow f\langle d
\rangle \leftrightarrow f\langle d' \rangle \leftrightarrow
f\langle c \rangle $ is also a triangle in $\htc(S)$. By the
observation above we have $f\langle c \rangle- f\langle d \rangle=
f\langle c \rangle-f\langle d' \rangle$.

Suppose now that $D'$ is an arbitrary simple closed curve on $S$
which is dual to $C$ and is disjoint from all $C_i$ for $i\geq 2$.
Then $d$ and $d'$ are two vertices of the complex $X_C$.
Since this complex is connected by Lemma
\ref{lemmaxv}, there is a sequence $d = d_1, d_2, \ldots,
d_n=d'$ of vertices in $X_C$ such that $d_i$ is connected to $d_{i+1}$ by
an edge for all $i=1,2,\ldots,n-1$. By the previous paragraph, we
have $f\langle c \rangle-f\langle d_i \rangle= f\langle c
\rangle-f\langle d_{i+1} \rangle$. It follows that $f\langle c
\rangle-f\langle d \rangle= f\langle c \rangle-f\langle d'
\rangle$.

This proves the lemma.
\end{proof}

\begin{lemma} \label{lemma2}
For a nonseparating simple closed curve $C$, the definition
of $\tilde{f}(c)$ is independent of all choices.
\end{lemma}

\begin{proof}
Suppose that $\{C_2,C_3,\ldots,C_g,D\}$ and $\{C'_2,C'_3,\ldots,C'_g,D'\}$
are two choices in the definition of $\tilde{f}(c)$. We must prove that
both choices give rise to the same result. More precisely, if
$v,w,v'$ and $w'$ denote the cut systems
$\langle c,c_2,\ldots,c_g \rangle $, $\langle d,c_2,\ldots,c_g \rangle $,
$\langle c,c'_2,\ldots,c'_g \rangle $ and $\langle d',c'_2,\ldots,c'_g \rangle $
respectively such that $v\leftrightarrow w$ and $v'\leftrightarrow w'$,
then we must show that $f(v')-f(w')=f(v)-f(w)$.

If $g=1$ then there are no $C_i$ and $C'_i$ and the conclusion of the lemma
follows from Lemma~\ref{lemma1}. So we assume that $g\ge 2$.

Suppose first that $v'$ is connected by an edge to $v$. Therefore
there are elements $c_{i_0}\in v$ and $c'_{j_0}\in v'$ such that
$C_{i_0}$ and $C'_{j_0}$ intersect transversely at one point and
$v-\{ c_{i_0}\} = v'- \{ c'_{j_0} \}$. After reindexing if
necessary we can assume that $c_{i_0}=c_2 $ and $c'_{j_0}=c'_2$,
so that  $c'_i=c_i$ for $i\geq 3$. Let $E$ be a simple closed
curve dual to $C$ and disjoint
from all $C_i$ and $C_2'$. Let $t = \langle e, c_2, c_3, \ldots,
c_g \rangle$ and $t' = \langle e, c'_2, c_3', c_4', \ldots, c_g'
\rangle = \langle e, c'_2, c_3, c_4, \ldots, c_g \rangle$. Since
$v\leftrightarrow v'\leftrightarrow t'\leftrightarrow t \leftrightarrow v$
form a rectangle in $\htc (S)$ and $f$ is
an automorphism,
$f(v)\leftrightarrow f(v')\leftrightarrow f(t')\leftrightarrow f(t) \leftrightarrow f(v)$
is a rectangle in $\htc (S)$. Then it is easy to see that $f(v) - f(t) =
f(v') - f(t')$. By using Lemma \ref{lemma1}, we obtain
$f(v) - f(t) =f(v) - f(w)$ and $f(v') - f(t')=f(v') - f(w')$.
Therefore, we get the desired result $f(v) - f(w)= f(v') - f(w')$.

\begin{figure}[htb]
\begin{center}
\epsfxsize=10cm
\epsfbox{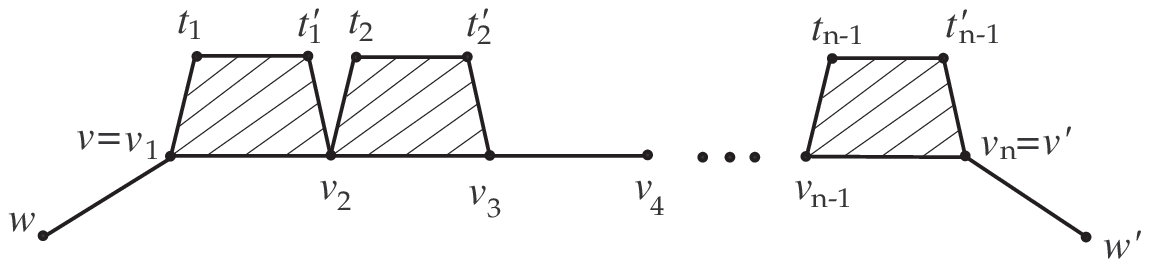}
\caption{A path in $\htc (S)$}
\label{ucgenler}
\end{center}
\end{figure}

Let us now consider the general case. Let $R$ denote the surface
obtained by cutting $S$ along the curve $C$. Thus $R$ is a surface
of positive genus. Since all $C_i$ and
$C'_i$, $i\geq 2$, are disjoint from $C$, we can consider them as
curves on $R$. Now $V=\langle   c_2,c_3, \ldots, c_g \rangle$ and
$V'=\langle c'_2,c'_3, \ldots, c'_g \rangle$ are two cut systems
on $R$. Since the Hatcher-Thurston complex $\htc (R)$ is connected
Theorem~\ref{thm:htw}, there is a sequence
$V=V_1,V_2,V_3,\ldots,V_n=V'$ of cut systems on $R$ such that
$V_i$ is connected by an edge to $V_{i+1}$. If we denote by $v_i$
the cut system on $S$ obtained from $V_i$ by adding $c$, we get a
path $v = v_1, v_2, \ldots, v_n = v'$ in $\htc (S)$. For each
$i=1,2,\ldots,g-1$, choose vertices $t_i$ and $t'_i$ as in the
previous paragraph such that $v_i\leftrightarrow
v_{i+1}\leftrightarrow t'_i\leftrightarrow t_i \leftrightarrow
v_i$ is a rectangle. We showed above that $f(v_i)- f(t_i) =
f(v_{i+1}) - f(t'_i)$. By Lemma~\ref{lemma1}, we also have
$f(v_{i+1}) - f(t'_i)= f(v_{i+1}) - f(t_{i+1})$. It follows that
$f(v_1)- f(t_1) = f(v_n) - f(t'_{n-1})$. Now the conclusion
$f(v)-f(w)=f(v')-f(w')$ follows from Lemma~\ref{lemma1}.

This completes the proof of the lemma.
\end{proof}

\begin{lemma}
\label{lemmaintone} Let $c, d$ be the isotopy classes of two
nonseparating simple closed curves $C$ and $D$
such that $i(c, d)= 1$. Then $i(\tilde{f}(c), \tilde{f}(d))= 1$.
\end{lemma}

\begin{proof}
It is easy to see that we can find nonseparating simple closed curves
$C_2, C_3,\ldots,C_g $ on $S$ such that
$v = \langle c, c_2,\ldots , c_g \rangle$ and $w = \langle d,
c_2,\ldots , c_g \rangle$ are two vertices in $\htc(S)$. Since
the geometric intersection of $c$ and $d$ is 1, we see that $v$
and $w$ are connected by an edge in $\htc(S)$. Since $f$ is
an automorphism $f(v)$ and $f(w)$ are also connected by an edge in
$\htc(S)$. From the definition of $\tilde{f}$,
we have $\{ \tilde{f}(c) \} = f(v) - f(w) $ and $\{ \tilde{f}(d) \} =f(w) - f(v) $.
Since $f(v)$ and $f(w)$ are connected by an edge, we
conclude that $i(\tilde{f}(c), \tilde{f}(d))=1$.\end{proof}

\begin{lemma}
\label{lemmatilde} If $f$ and $h$ are two automorphisms of
$\htc(S)$ and if $c$ is the isotopy class of a nonseparating simple closed curve $C$,
then $\widetilde{fh}(c)=\tilde{f} (\tilde{h}(c))$.
\end{lemma}
\begin{proof}
Let us choose vertices
$c_2,c_3,\ldots,c_g$ and $d$ in $\mathcal{G}(S)$ such that
$v=\langle c,c_2,\ldots, c_g\rangle$ and
$w=\langle d,c_2,\ldots, c_g\rangle$ are distinct vertices in the complex
$\htc(S)$ which are connected by an edge; $v\leftrightarrow w$.
Then $\{ \tilde{h}(c) \} =h(v)-h(w)$ and $\{ \widetilde{fh}(c) \} =fh(v)-fh(w)$.
Since $h(v)\leftrightarrow h(w)$, we can use
these vertices to define $\tilde {f} (\tilde{h}(c))$:
\begin{eqnarray*}
\{ \, \tilde {f} (\tilde{h}(c))\, \} &=& f(h(v))-f(h(w))\\
&=& (fh)(v)-(fh)(w)\\
&=& \{\, (\widetilde  {fh})(c)\,\}.
\end{eqnarray*}
\end{proof}

\begin{proposition}
\label{prop1}
The mapping $\tilde{f}$ is an automorphism of the
graph $\mathcal{G}(S)$.
\end{proposition}

\begin{proof} For a vertex $c$ in $\mathcal{G}(S)$, $\tilde{f}(c)$ is
well-defined. Therefore we have a well-defined  map $\tilde{f} : \mathcal{G}(S) \to \mathcal{G}(S)$.
If two vertices $c,d$ are connected by an edge in
$\mathcal{G}(S)$, then $i(c,d)=1$. By Lemma~\ref{lemmaintone},
$i(\tilde{f}(c), \tilde{f}(d))= 1$. Therefore, $\tilde{f}$ is simplicial.

Let $h\in {\rm Aut}\,\htc (S)$ be the inverse of $f$. Then
$\tilde{f}\tilde{h}$ and $\tilde{h}\tilde{f}$ are both the
identity automorphisms, because it can be shown that if $I\in {\rm Aut}\,\htc (S)$
denote the identity, then $\tilde{I}(c)=c$ for all nonseparating simple closed curve $C$.
We conclude that $\tilde{f}: \mathcal{G}(S) \to \mathcal{G}(S)$ is a bijection.\end{proof}

In the following proposition we will prove that $\tilde{f}$
preserves geometric intersection zero, and hence also is an
automorphism of $\mathcal{N}(S)$ for closed surfaces.

\begin{proposition}
\label{prop2} If $S$ is a closed surface of genus at least two,
then the mapping $\tilde{f} : \mathcal{N}(S)\to \mathcal{N}(S)$
is an automorphism.
\end{proposition}

\begin{proof} By the previous proposition, $\tilde{f}$ is a bijection.
So, it is enough to show that
$\tilde{f}$ is a simplicial map on $\mathcal{N}(S)$. Let $a, b$ be
two distinct vertices of $\mathcal{N}(S)$, which have disjoint
representatives $A$ and $B$ on $S$ respectively. We will consider
the following two cases:

Case i: If $S_{A \cup B}$ is connected, then $\{a, b\}$ can be
completed to a vertex $v$ in $\htc(S)$. Then, since $f(v)$ is a
vertex in $\htc(S)$ and $\tilde{f}(a), \tilde{f}(b) \in f(v)$, we
see that $\tilde{f}(a)$ and $\tilde{f}(b)$ have disjoint
representatives on $S$.

Case ii: If $S_{A \cup B}$ is not connected, then we complete $A$
and $B$ to a curve configuration as shown in Figure 3, by taking a
maximal chain $\{C_1, ..., C_{2g+1}\}$ with $i(c_i, c_{i+1})=1$,
$i(c_i, c_j)=0$ for $|i-j| > 1$, $c_i \in \mathcal{N}(S)$ as shown
in the figure for $g=4$ case (similar chains can be chosen in the
other cases). Notice that $S_{C_i \cup C_j}$ is connected for any
$i, j$. So, if $i(c_i, c_j)$= 0 then $i(\tilde{f}(c_i),
\tilde{f}(c_j))$= 0 by the first case. If $i(c_i, c_j)$= 1, then
$i(\tilde{f}(c_i), \tilde{f}(c_j))$= 1 by Lemma \ref{lemmaintone}.
Hence $\{\tilde{f}(c_1), ..., \tilde{f}(c_{2g+1})\}$ is a maximal
chain on $S$.

\begin{figure}[htb]
\begin{center}
\epsfxsize=2.7in \epsfbox{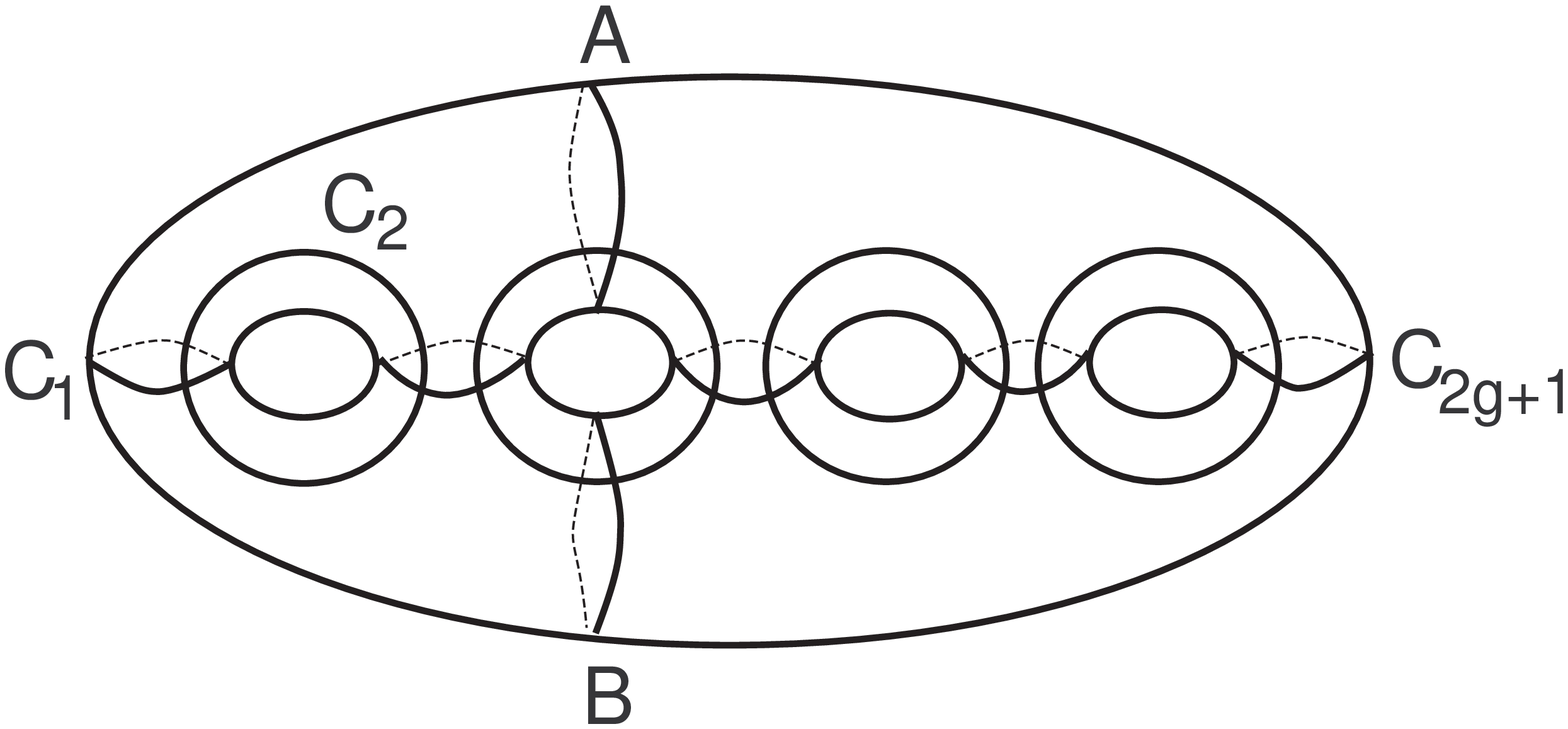} \caption{A, B and a chain}
\end{center}
\end{figure}

We also have that $S_{A \cup C_i}$ is connected for any $i$. So,
if $i(a, c_i)$= 0 then $i(\tilde{f}(a), \tilde{f}(c_i))$= 0 by the
first case. If $i(a, c_i)$= 1, then $i(\tilde{f}(a),
\tilde{f}(c_i))$= 1 by Lemma \ref{lemmaintone}. Similarly since
$S_{B \cup C_i}$ is connected for any $i$, if $i(b, c_i)$= 0 then
$i(\tilde{f}(b), \tilde{f}(c_i))$= 0 by the first case. If $i(b,
c_i)$= 1, then $i(\tilde{f}(b), \tilde{f}(c_i))$= 1 by Lemma
\ref{lemmaintone}.

Note that $i(a,c_{2k})=i(b,c_{2k})=1$ for some integer $k\in \{
2,3,\ldots ,g-1 \}$ and the intersection numbers of $a$ and $b$
with any other $c_i$ is $0$. Therefore,
$i(\tilde{f}(a),\tilde{f}(c_{2k}))
=i(\tilde{f}(b),\tilde{f}(c_{2k}))=1$ and the intersection numbers
of $\tilde{f}(a)$ and $\tilde{f}(b)$ with any other
$\tilde{f}(c_i)$ is $0$.

Let $C_i' \in \tilde{f}(c_i)$, $A' \in \tilde{f}(a)$ and $B' \in
\tilde{f}(b)$ such that all the curves $C_i'$, $A'$ and $B'$
intersect minimally with each other for each $i$. Since $A$ and
$B$ are dual to $C_{2k}$, $A'$ and $B'$ are dual to $C'_{2k}$ by
Lemma \ref{lemmaintone}. Since curves $A'$ and $B'$ are disjoint
from the chains $C'_1\cup C'_2\cup \cdots \cup C'_{2k-1}$ and
$C'_{2k+1}\cup C'_{2k+2}\cup \cdots \cup C'_{2g+1}$ and since the
complement of these two chains is the union of two annuli, the
(distinct) curves $A'$ and $B'$ must be disjoint. Because up to
isotopy there are only two simple closed curves on the disjoint
union of two annuli and they are disjoint. Therefore,
$i(\tilde{f}(a), \tilde{f}(b)) = 0$. This shows that $\tilde{f}$
is a simplicial map on $\mathcal{N}(S)$. Since $\tilde{f}$ is 1-1
and onto, it is an automorphism of $\mathcal{N}(S)$.\end{proof}

\noindent {\bf Remark:} If $S$ is a closed surface of genus at
least two, by using Proposition \ref{prop2} and the results in
\cite{Ir3}, we see that $f$ is induced by a homeomorphism of $S$.

\section{Automorphisms of $\htc(S) $ and mapping class group}
\label{section4}

In this final section, we state and prove the  main result. We then give a corollary to
the main theorem and comment on other possible but similar proofs of the main theorem.

\begin{theorem}
\label{thmmain}
Let $S$ be a compact, connected, orientable surface of genus at
least one. Then the mapping $\varphi: {\rm Aut\,}\htc(S)\rightarrow  {\rm Aut\,}\mathcal{G}(S)$
given by $f\mapsto \tilde{f}$ is an isomorphism.
\end{theorem}

\begin{proof} By the results of the previous section,
$\varphi (f)=\tilde{f}$ is a well-defined automorphism of
Aut\,$\mathcal{G}(S)$, and Lemma~\ref{lemmatilde} shows that
$\varphi$ is a group homomorphism.

For an element $f\in {\rm Aut\,}\htc(S)$ if $\tilde{f}$ is the
identity automorphism of $\mathcal{G}(S)$, then it follows from
$f(\langle c_1,c_2,\ldots, c_g \rangle) = \langle
\tilde{f}(c_1),\tilde{f}(c_2),\ldots,\tilde{f}(c_g) \rangle $ that
$f$ acts trivially on $\htc(S)$. Hence, $\varphi$ is one-to-one.

If $h$ is an automorphism of $\mathcal{G}(S)$, then $h$ is induced
by a homeomorphism $F$ of the surface $S$ by the results given in
\cite{Sc}. Now $F$ induces an automorphism $f$ of $\htc(S)$ and
$\tilde{f}=h$. Hence, $\varphi$ is an isomorphism. This completes
the proof of the theorem.
\end{proof}

\begin{corollary}
Let $S$ be a compact, connected, orientable surface of genus
$g\geq 1$ with $r\ge 0$ boundary components. If $(g,r)\neq
(1,0),(1,1),(1,2),(2,0)$, then we have ${\rm Aut\,} \htc(S) \cong
\mod_S^*$. If $(g,r)$ is one of $(1,0),(1,1),(1,2),(2,0)$, then we
have ${\rm Aut\,} \htc(S) \cong \mod_S^*/\Z_2$. That is, in all
cases ${\rm Aut\,}\htc(S) \cong \mod_S ^* /\mathcal{C}(\mod_S
^*)$.
\end{corollary}
\begin{proof}
The proof follows from Theorem~\ref{thmmain}
and Theorem~\ref{thmsc}.
\end{proof}

This gives us another evidence to the following conjectural
statement; the automorphism group of any natural complex on curves
is isomorphic to the extended mapping class group in generic
cases. By using our results given in this paper, and the main
results of \cite{Ir1}, \cite{Ir2}, \cite{Ir3}, \cite{Iv},
\cite{Sc}, \cite{M}, we see that for most of the compact,
connected, orientable surfaces, we have ${\rm Aut\,} \htc(S) \cong
{\rm Aut\,} \mathcal{N}(S) \cong {\rm Aut\,} \mathcal{C}(S) \cong
{\rm Aut\,} \mathcal{G}(S) \cong {\rm Aut\,} \mathcal{P}(S)$ where
$\mathcal{P}(S)$ is the pants complex.\\

If $S$ is a closed surface of genus at least two and $f$ is an
automorphism of $\htc(S)$, by using the techniques given in this
paper, in particular Lemma \ref{lemmaintone}, we can see that $f$
induces an automorphism $\tilde{f}_*$ on $\mathcal{C}(S)$ by
extending $\tilde{f}$ over the nontrivial separating curves on $S$
by using chains on two subsurfaces that the separating curves
separate. We do the following: Let $C$ be a nontrivial separating
curve on $S$. Since $g \geq 2$, $C$ separates $S$ into two
subsurfaces $S_1, S_2$, and both of $S_1, S_2$ have genus at least
one. We take a chain on $S_1$, $\{a_1,a_2, \ldots , a_m\}$ with
$i(a_{i}, a_{i+1})=1$, $i(a_{i}, a_{j})=0$ for $|i-j| > 1$, $a_i
\in \mathcal{N}(S)$, such that $S_1 \cup \{c\}$ is a regular
neighborhood of $A_1 \cup A_2 \cup \cdots \cup A_n$ where $A_i\in
a_i $ and $A_i$'s intersect minimally. Since $\tilde{f}$ preserves
disjointness and intersection one property, we can see that the
chain $\{a_1,a_2,\ldots , a_m \}$ is mapped by $\tilde{f}$ into a
similar chain, $\{\tilde{f}(a_1), \ldots , \tilde{f}(a_m) \}$ with
$i( \tilde{f}(a_{i}), \tilde{f}(a_{i+1}))=1$, $i(\tilde{f}(a_{i}),
\tilde{f}(a_{j}))=0$ for $|i-j| > 1$. Let $A_i' \in
\tilde{f}(a_i)$ such that any two elements in $\{A_1', \ldots ,
A_m'\}$ have minimal intersection with each other. Let $M$ be a
regular neighborhood of $A_1' \cup A_2' \cup \ldots \cup A_m'$.
Then it is easy to see that $M$ is homeomorphic to $R_1 \cup c$.
Let $C'$ be the boundary of $M$. We define $\tilde{f}_* (c) =
[C']$ (See \cite{Ir3} for well definedness). This extends
$\tilde{f}$ to a simplicial map $\tilde{f}_*$ on $\mathcal{C}(S)$.
It can be shown that $\tilde{f}_*$ is an automorphism on
$\mathcal{C}(S)$.\\

\newpage
{\bf Acknowledgments}\\

We would like to thank Joan Birman, Nikolai Ivanov, John McCarthy
and Dan Margalit for their interest in this work and for valuable
comments about this paper.

\end{document}